\documentstyle [11pt,amssymb]{article}

\begin {document}

\title {On Markovian Cocycle Perturbations in Classical
and Quantum Probability}

\author {G.G. Amosov\thanks {The work is partially supported by
INTAS-00-738}}

\maketitle

\abstract {We introduce Markovian cocycle perturbations
of the groups of transformations associated with
the classical and quantum stochastic processes with
stationary increments, which are characterized by
a localization of the perturbation to the algebra
of events of the past. It is namely the definition one
needs because the Markovian perturbations of the
Kolmogorov flows associated with the
classical and quantum noises result in the perturbed
group of transformations which can be decomposed in the sum of
a part associated with
deterministic stochastic processes lying
in the past and a part associated with
the noise isomorphic to the initial
one. This decomposition allows to obtain some analog
of the Wold decomposition for classical stationary
processes excluding a nondeterministic part of the process
in the case of the stationary quantum stochastic
processes on the von Neumann factors which are
the Markovian perturbations of the quantum noises.
For the classical stochastic process with noncorrelated
increaments it is constructed the model of Markovian
perturbations describing all Markovian cocycles up to
a unitary equivalence of the perturbations. Using this
model we construct Markovian cocyclies transformating
the Gaussian state $\rho $ to the Gaussian states
equivalent to $\rho $.
}

\section*{Introduction}

It is well-known that every stochastic process with finite
second moments and continuous in square mean can be considered
as a continuous curve in the Hilbert space.
In this framework properties of the process such as
stationarity and noncorrelativity of the increments
appear as an invariance under the action of the group of unitaries
in the Hilbert space and an orthogonality of the curve increments
correspondingly. Thus to investigate the stochastc process
one can use the functional analysis techniques.
This approach was introduced by A.N. Kolmogorov who
considered in \cite {Kol1, Kol2}
a classification problem for the equivalence classes of
continuous curves $\xi =(\xi _t)_{t\in \mathbb R}$ in
a Hilbert space $H$, which are invariant with respect to
a strong continuous one-parameter group of unitaries
$U=(U_t)_{t\in \mathbb R}$, with the transformations of
equivalence defined by the formula
$$
\tilde \xi _t=W\xi _t+\eta ,
\eqno (*)
$$
where $W$ and $\eta $ are a unitary operator in $H$ and an
element of $H$ correspondingly. Indeed, the continuous curve
$\tilde \xi = (\tilde \xi _t)_{t\in \mathbb R}$ is invariant
with respect to the group of unitaries $\tilde U=(W_tU_t)_{t\in
\mathbb R}$, where the one-parameter family of unitaries
$$
W_t=WU_tW^*U_{-t},\ t\in \mathbb R,
\eqno (**)
$$
satisfies the condition of multiplicative {\it $U$-cocycle},
$$
W_{t+s}=W_tU_tW_sU_{-t},\ s,t\in \mathbb R.
\eqno (***)
$$
The multiplicative cocycle generated by a unitary operator
$W$ as in $(**)$ is said to be {\it $U$-
coboundary}. Notice that not every cocycle $(***)$ is
coboundary $(**)$. Here we use ordinary
definitions of the cohomologies of groups theory for $1$-cocycle
and $1$-coboundary associated with the standard
bar-resolvent of the group $\mathbb R$ with values in
the multiplicative group ${\cal U}(H)$ of all unitaries
in the Hilbert space $H$ with the module structure
defined by the group action $x\to U_txU_{-t},\
x\in {\cal U}(H),\ t\in {\mathbb R}$ (see, f.e.,
\cite {Bra, Gui}).

Suppose that there exists a continuous curve
$\xi =(\xi _t)_{t\in {\mathbb R}}$ which is invariant
with respect to the group $U$ and the
increments of $\xi $ are orthogonal such that
$\xi _{t_1}-\xi _{s_1}
\perp \xi _{t_2}-\xi _{s_2}$ for all disjoint intervals
$(s_1,t_1)\cap (s_2,t_2)=\emptyset $. The curves of
such type was called in \cite {Kol2} by the Wiener spirales.
Every multiplicative $U$-cocycle $(W_t)_{t\in {\mathbb
R}}$, which is strong continuous in $t$, defines a new
strong continuous one-parameter group of unitaries
$\tilde U=(W_tU_t)_{t\in {\mathbb R}}$.
We shall call $\tilde U$
{\it a cocycle perturbation} of the group $U$. If
$(W_t)_{t\in {\mathbb R}}$ is a coboundary, then
for the cocycle perturbation $\tilde U$ there also
exists the invariant curve $\tilde \xi $ which is
the Wiener spiral. Indeed, this curve can be constructed by
the formula $(*)$, where one must substitute for
$W$ the unitary operator generating the coboundary
by means of $(**)$.
For arbitrary cocycle
$(W_t)_{t\in {\mathbb R}}$ which is not a coboundary
it is possible that there doesn't exist the Wiener
spiral which is ivariant with respect to the cocycle perturbation
$\tilde U$. In this paper we introduce a subset $(M)$ of
the set of all cocycles such that given a cocycle from $(M)$
generates the cocycle perturbation possessing the invariant
Wiener spiral if the perturbing group satisfied this condition.
The class of cocycle perturbations we propose is important
because in applications (see \cite {Par, Hol0, Hol1})
it is often posed the problem of representing the quantum
stochastic process as a projection of the cocycle perturbation
of the quantum noise which is the Wiener spiral from the
viewpoint of functional analysis.
Using in this framework cocycles satisfying
the quantum stochastic differential equations belong
to the subset $(M)$ we introduce.
Cocycles from $(M)$ we shall call Markovian.
Using of the term "Markovian cocycle" follows from
\cite {Acc1, Acc2}.  The Markovian cocycles of
\cite {Acc1, Acc2} are connected with the classical
Markov property in the sense that the perturbation of
the Markovian stochastic process by the
Markovian cocycle is also the Markovian stochastic
process with the same algebras of the present, the future
and the past. Our definition distinguishes from
the definition in \cite {Acc1, Acc2}, nevertheless,
under an appropriate interpretation
the perturbations we consider are Markovian in the sense of
the definition in the cited paper descrbing more wide
class of perturbations. We consider it as
a sufficiant motivation to use this terminology.

The model example of the group of unitaries $U$
possessing the invariant Wiener spiral is
the group of shifts in the Hilbert space
$H=L^2({\mathbb R})$ defined by the formula
$(U_t\eta )(x)=\eta (x+t),\
\eta \in H$. In fact, the set of functions
$\xi _t(x)=1,\ x\in [-t,0],\ \xi _t(x)=0,\
x\notin [-t,0],$ is the Wiener spiral which is
invariant with respect to $U$ (see \cite {Kol2}).
We shall call by "a past" and "a future" of the system
the subspaces $H_{t]}$ and $H_{[t}$ containing
functions with the support belonging to $[t,+\infty )$ and
$(-\infty ,t]$ correspondingly.
Notice that
$H_{t]}=U_{t}H_{0]},\  H_{[t}=U_{t}H_{[0},\
t\in {\mathbb R}$. In the case we are considering
we call the cocycle $W$ {\it Markovian}
if the restriction of the unitary operator
$W_t$ to the subspace of the future
$H_{[t}$ is an identity transformation for
every fixed $t\geq 0$. Particularily, this definition
garantees that the subspace of the past $H_{0]}$
is invariant with respect to $W_{-t}$ for all
$t\geq 0$, which allows to consider the restriction
of the cocycle perturbation
$\tilde U_{-t}|_{H_{0]}},\ t\geq 0$.
Let $\tilde U$ be a cocycle perturbation of
$U$ by the cocycle which is Markovian in the
sense of our definition. We shall show that
every such perturbation can be represented in
the form
$\tilde U=\tilde U^{(1)}\oplus
\tilde U^{(2)}$, where $\tilde U^{(1)}$
is arbitrary group of unitaries in the
subspace of the past $H_{0]}$
and $\tilde U^{(2)}$ is unitary equivalent
to the initial group of shifts.
This representation can be named the Wold
decomposition of the cocycle perturbation
of $\tilde U$. The groups $\tilde U^{(1)}$
and $\tilde U^{(2)}$
can be interpreted as associated with
deterministic and nondeterministic
parts of the process.

Under a quantum stochastic process
we mean (see \cite {Par, Hol0, Hol1})
the strong continuous one-parameter family
$x=(x_t)_{t\in \mathbb R}$ of linear (nonbounded in
general) operators in a Hilbert space.
Accordingly to this definition stationary quantum stochastic process
$x=(x_t)_{t\in {\mathbb R}}$ can be defined by the condition
$x_t=\alpha _t(x^0)+x^1,$ $t\in {\mathbb R},$ where $\alpha
$ and $x^0,x^1$ are certain group of automorphisms and
two fixed linear operators correspondingly.
In this way, a quantum stochastic process with stationary increments
is the continuous operator-valued curve which is invariant with
respect to certain group of automorphisms $\alpha $.  In the quantum
probability theory a role of the $\sigma $-algebras of events
associated with the stochastic process is played the von Neumann
algebras generated by increments of the process, i.e. the ultra-weak
closed algebras of bounded operators in a Hilbert space
with the commutant (the set of bounded operators
commuting with operators of the algebra) defined by
the condition of commuting with increments of the process
$x$ (see, f.e. \cite {Hol0, Hol1, Brat}).
Every classical stochastic process $\xi
$ consisting of the random variables $\xi _t\in L^{\infty }(\Omega )$
can be considered as quantum, where the operators $x_t=M_{\xi _t}$
forming the process are the operators of multiplications by
the functions $\xi _t$ in the Hilbert space $L^2(\Omega )$.

Remember that the flow
$\{T_t,\ t\in {\mathbb R}\}$ on the probability space
$(\Omega ,{\cal M},\mu )$ is said to be a Kolmogorov flow
(see \cite {Kol3}) if
there exists a $\sigma $-algebra of events ${\cal M}_{0]}\subset
{\cal M}$ such that $T_t{\cal M}_{0]}={\cal M}_{t]},$
$$
{\cal M}_{s]}\subset {\cal
M}_{t]},\ s<t,
\eqno (K_1)
$$
$$ \cup _t{\cal M}_{t]}={\cal M},
\eqno (K_2)
$$
$$
\cap _t{\cal M}_{t]}=\{\emptyset ,\Omega \}.
\eqno (K_3)
$$
In \cite {Emc} a notion of the Kolmogorov flow was
transfered into
quantum probability. Under this construction
the $\sigma $-algebras of events
${\cal M}_{t]}$
are replaced by the corresponding von Neumann algebras.
Notice that in the well-known monograph
\cite {Ibr}, where it was investigated the conditions
on the spectral function of process, which result in
the Kolmogorov flow, it is used a term "completely
nondeterministic process" for the stochastic process
generating the Kolmogorov flow, while a term "Kolmogorov
flow" is not used anywhere.
Every classical or quantum process with independent
(in classical sense) increments results in the Kolmogorov
flow. We define Markovian cocycle perturbations of
the classical and quantum Kolmogorov flows such that
the perturbed flows contain the parts which are isomorphic
to the initial Kolmogorov flow. Notice that in the model
situation of the Hilbert space
$L^2({\mathbb R})$ we considered above the Kolmogorov
flow can be associated with the flow of shifts.
Thus the possibility to exclude in the perturbed dynamics
the part being the Kolmogorov flow, one can consider as
some analogue of the Wold decomposition allowing to exclude
a nondeterministic part of the process.

For every classical or quantum stochastic process a
great role is played the set of all (not nessesarily linear)
functionals of the process. Particularily, for the Wiener
process it is defined the Wiener-Ito decomposition of
the space of all $L^2$-functionals in the orthogonal sum,
which allows to solve effectively the stochastic differential
equations (see \cite {Hid}). Notice that from the viewpoint of
the theory of the cohomologies of groups, the stochastic
process with stationary increments determining the continuous
curve being invariant with respect to the group of
transformations $U$ is an additive $1-U$-cocycle.
It is naturaly to consider the ring of cohomologies of
all degrees generated by the $1-U$-cocycle of such type
which can be interpreted as the space of all (nonlinear)
functionals of the initial stochastic process.
We show that the Markovian cocycle perturbations
we introduced define homomorphisms of this ring
of cohomologies.

This paper is organized as follows.
In Part 1 we define a class of the Markovian
$U$-cocycles such that the group of unitaries
$\tilde U$ which is obtained through the perturbation
by the cocycle of such type determines a continuous
curve $\tilde \xi $ being invariant with respect to
$\tilde U$ and connected with the initial continuous
curve $\xi$ by the formula $(*)$, where $W$ is (nonunitary
in general) isometrical operator satisfying the additional
property of the localization of action to "the past".
Such isometrical operators we call Markovian and prove
that every Markovian operator is associated with
certain perturbation by a Markovian cocycle.
The investigation in details we give in Part 2 shows
that the Markovian cocycle we introduced determines
the Wold decomposition for the cocycle perturbation.
In Part 3 we construct a model of the Markovian cocycle
for the stochastic process with independent increments.
The model we gave allows to construct the Markovian
cocycles with the property $W_t-I\in s_2$
(the Hilbert-Schmidt class) which translate
the fixed Gaussian measure to the equivalent Gaussian measures.
In Part 4 it is given the basic notion on the theory of
Kolmogorov flows in classical and quantum probability.
In Part 5 we define the rings of cohomologies generated
by additive $1$-cocycles and show in examples of
the Wiener process and the quantum noise that the set of
all functionals of the stochastic process with stationary
increments can be considered as a ring of the cohomologies
of the group which is composed by shifts of the increments in time.
Moreover, we define a Markovian perturbation of the group
resulting in a homomorphism of the ring of cohomologies of
the group. Notice that for quantum stochastic processes
$x=(\xi _t=x_t)_{t\in \mathbb
R}$ it is also possible to define transformations of
the form $(*)$, where one must take a morphism for
$W$ and a linear operator for $\eta $.
In Part 6 we introduce the Markovian cocycle perturbations
of the quantum noises being a generalization of
the classical processes with independent stationary
increments for the quantum case, which result in
transformations of the quantum stochastic processes
of the form $(*)$.
The Markovian perturbations we introduce determine
homomorphisms of the ring of cohomologies associated
with the stochastic proocess in the sense of Part 5.
The techniques involved in Part 5 allows to obtain
in Part 6 some analogue of the Wold decomposition
for the classical stochastic processes in the quantum
case, which permits to exclude a nondeterministic part
of the process.

\section {Stochastic processes with stationary
increments as curves in a Hilbert space,
which are invariant with respect to the group
of transformations
}

Let $\xi =(\xi_t)_{t\in \mathbb R}$ be a continuous in
square mean stochastic process with stationary increments
on the probability space $(\Omega
,\Sigma,\mu)$. Without lost of a generality
we can suppose that the condition $\xi (0)=0$ holds. Then in
the Hilbert space $H^{\xi}$ generated by the increments $\xi
_t-\xi _s,\ s,t\in \mathbb R,$ it is defined a strong
continuous group of unitaries $U=(U_t)_{t\in \mathbb R}$
shifting the increments in time such that
$U_t(\xi _s-\xi _r)=\xi
_{s+t}-\xi _{r+t},\ s,t,r\in \mathbb R$, where $\xi $ satisfies
the condition of (additive) $1-U$-cocycle, i.e. $\xi _{t+s}=\xi
_t+U_t\xi _s,\ s,t\in \mathbb R$.  Denote $H$ the Hilbert space
with the inner product $(\xi,\eta)=\mathbb E (\xi
\overline \eta )$ generated by the classes of equivalency of
the random variables $\xi ,\eta $ in the space $(\Omega ,\mu)$,
which possess finite second moments, with repect to the
Hilbert norm associated with the expectation  $\mathbb E$.
We shall identify the random variables with the elements of the
Hilbert space $H$.  A stochastic process $\xi $ can be considered as
a curve in the Hilbert space $H^{\xi }$, which is invariant
with repect to the action of the group $U$ (see \cite {Kol1},\cite
{Kol2}).  The space $H^{\xi }$ is a subspace of $H$ but doesn't
coincide with it in general.
Let $\tilde U$ be arbitrary continuation of
$U$ to a strong continuous group of unitaries in the
Hilbert space $H$. Then $\xi $ can be considered also as
the curve in $H$, which is invariant with respect to
the action of the group
$\tilde U$. Denote ${\cal I}(\xi )$ the set containing
all possible strong continuous groups of unitaries in
$H$ such that the stochastic process $\xi $ is invariant
with respect to them.
One can define in the space $H^{\xi }$ an increasing family of
subspaces $H_{t]}^{\xi }$
generated by the increments $\xi _s-\xi_r,\ s,r\leq t,$ associated
with "a past" before the moment $t$, and a decreasing family of
subspaces $H_{[t}^{\xi }$ generated by the increments $\xi _s-\xi
_r,\ s,r\geq t$ associated with "a future" after the moment $t$
such that $H^{\xi }=\vee _tH_{t]}^{\xi }=\vee _tH_{[t}^{\xi }$.
Notice that for the processes with noncorrelated increments
the subspaces
$H_{t]}^{\xi }$ and $H_{[t}^{\xi }$ are orthogonal.
Fix the group
$U\in {\cal I}(\xi )$. The strong continuous one-parameter family of
unitaries $W=(W_t)_{t\in \mathbb R}$ in $H$ is said to be {\it a
(multiplicative) $U$-cocycle} if the following condition holds, $$
W_{t+s}=W_tU_tW_sU_t^*,\ s,t\in {\mathbb R}, W_0=I.
\eqno (C)
$$
The cocycle $W$ is called {\it Markovian} under the condition
$$
W_tf=f,\ f\in H_{[t}^{\xi },\ t>0.
\eqno (M)
$$
The property $(C)$
exactly means that the strong continuous one-parameter family of
unitaries $\tilde U=(W_tU_t)_{t\in \mathbb R}$ forms
a group.
We consider Markovianity as
a localization of the action of the cocycle $W$ to the subspace
of the past. Moreover the Markovian property $(M)$
preserves "a causality" such that "the future" of the system
is not disturbed. Notice that our definition of
a Markovian cocycle is based on the analogous definition
introduced in \cite {Acc1, Acc2} in considerably more general
case. We defer the examples of Markovian cocycles to Part 3,
where it is given a model of the Markovian cocycle for the group
of shifts on the line, which describes all cocycles up to
the unitary equivalence of perturbations.
Using $(C)$ we get
$I=W_{-t+t}=W_{-t}U_{-t}W_tU_{-t}^*,\ t>0$, such that
$W_{-t}=U_t^*W_tU_t,\ t>0$. Thus one can rewrite
$(M)$ in the form
$$
W_{-t}f=f,\ f\in H_{[0}^{\xi },\ t>0.
\eqno (M_1)
$$
Consider the stochastic process $\xi ^{(1)},\ \xi ^{(1)}(0)=0,$
being a continuous curve in the Hilbert space
$H$ and take $U^{(1)}\in {\cal I}(\xi ^{(1)})$.
Suppose $W=(W_t)_{t\in \mathbb R}$ is a multiplicative
$U^{(1)}$-cocycle in the space $H$.

{\bf Proposition 1.1.}{\it Let $W$ satisfy
the Markovian property $(M)$. Then
continuous in $t$ family
$\xi ^{(2)}_t=W_t\xi ^{(1)}_t,$ $t\leq 0,\
\xi ^{(2)}_t=\xi ^{(1)}_t,$ $t>0,$
is the stochastic process with stationary increments
which is a curve in $H$ being invariant with respect to
the group of unitaries
$U^{(2)}_t=W_tU^{(1)}_t,$ $t\in \mathbb R$.}

Proof.

Check that $\xi ^{(2)}_{t+s}=\xi ^{(2)}_t+
U^{(2)}_t\xi ^{(2)}_s,\ s,t\in \mathbb R$. In fact,
for $s,t\leq 0$ we obtain
$\xi ^{(2)}_{t+s}=W_{t+s}\xi ^{(1)}_{t+s}=
W_{t+s}\xi ^{(1)}_t+W_{t+s}U^{(1)}_t\xi ^{(1)}_s=
W_tU^{(1)}_tW_sU^{(1)}_{-t}\xi ^{(1)}_t+
W_tU^{(1)}_tW_s\xi ^{(1)}_s=
-W_tU^{(1)}_tW_s\xi ^{(1)}_{-t}+
U^{(2)}_t\xi ^{(2)}_s)=
-W_tU^{(1)}_t\xi ^{(1)}_{-t}+U^{(2)}_t\xi ^{(2)}_t=
\xi ^{(2)}_t+U^{(2)}_t\xi ^{(2)}_s$,
where we used the identity
$W_s\xi ^{(1)}_{-t}=\xi ^{(1)}_{-t},\ s,t\leq 0$,
which is correct due to $(M_1)$.
For $s,t>0$ the property we prove is true
because $\xi _t^{(2)}=\xi _t^{(1)},\ t\geq 0$,
by the definition. Notice that $\xi _t^{(2)}=\xi _t^{(1)}=
-U_t^{(1)}\xi _{-t}^{(1)}=-W_tU_t^{(1)}$ $W_{-t}U_{-t}^{(1)}
U_t^{(1)}\xi _{-t}^{(1)}=-W_tU_t^{(1)}W_{-t}\xi _{-t}^{(1)}=
-U_t^{(2)}\xi _{-t}^{(2)},\ t>0,$ by means of the cocycle
conditions $(C)$ for $W$. $\Box$

We shall call an isometrical operator
$R:H^{\xi ^{(1)}}\to H^{\xi ^{(2)}}$ Markovian if
$Rf=f,\ f\in H^{\xi ^{(1)}}_{[0}$.

{\bf Proposition 1.2.} {\it Given a Markovian cocycle
$W$ there exists the limit $\lim \limits _{t\to +\infty }
W_{-t}\eta =W_{-\infty }\eta ,\ \eta \in H^{\xi ^{(1)}}$,
such that
$W_{-\infty }$ is a Markovian isometrical operator
with the property
$W_{-\infty }f=W_{-s}f,\ f\in H^{\xi ^{(1)}}_{[-s},\ s\geq
0$.}

Proof.

Notice that $W_{-t-s}f=W_{-s}U^{(1)}_{-s}W_{-t}U^{(1)}_sf=
W_{-s}f,\ f\in H^{\xi ^{(1)}}_{[-s},$ $s,t\geq 0$, due to
the Markovian property in the form $(M_1)$.
Hence, the limit exists for the set of elements
$f\in H^{\xi ^{(1)}}_{[-s},\ s\geq 0$,
which is dense in $H^{\xi ^{(1)}}$.
Thus, the strong limit exists by the Banach-Steinhaus
theorem.
The limiting operator $W_{-\infty }$ is Markovian because
all operators $W_{-t},\ t\geq 0,$ satisfy this condition.
$\Box$

{\bf Proposition 1.3.} {\it In the space $H$ there exists
a Markovian isometrical operator $R$ with the property
$\xi ^{(2)}_t=R\xi ^{(1)}_t,\
t\in \mathbb R,$ if and only if
the stochastic processes $\xi ^{(1)}$ and $\xi ^{(2)}$ are
connected by a Markovian cocycle
$W=(W_t)
_{t\in \mathbb R}$ by the formula $\xi ^{(2)}_t=W_t\xi ^{(1)}_t,\
t\leq 0$.}

Proof.

Necessity.

Suppose that there exists a Markovian isometrical
operator $R$ such that
$\xi ^{(2)}_t=R\xi ^{(1)}_t,\ t\in \mathbb R$.
Check the condition $(M)$ for $W$,
$W_t(\xi ^{(1)}_r-\xi ^{(1)}_s)=U^{(2)}_tU^{(1)}_{-t}
(\xi ^{(1)}_r-\xi ^{(1)}_s)=U^{(2)}_t
(\xi ^{(1)}_{r-t}-\xi ^{(1)}_{s-t})=U^{(2)}_t
(\xi ^{(2)}_{r-t}-\xi ^{(2)}_{s-t})=
\xi ^{(2)}_r-\xi ^{(2)}_s,\
r,s\geq t>0$.
Here we used the Markovian property for
the operator $R$, which implies
$\xi ^{(2)}_t=R\xi ^{(1)}_t=\xi ^{(1)}_t,\ t\geq 0$.

Sufficiency.

Suppose that $\xi ^{(2)}_t=W_t\xi ^{(1)}_t,\ t\leq 0,$
and a cocycle $W$ is Markovian. Then,
due to Proposition 1.2, one can define a Markovian
isometrical operator $W_{-\infty }=s-\lim \limits
_{t\to +\infty }W_{-t}$.
Consider the process $\tilde \xi ^{(2)}_t=W_{-\infty }\xi _t^{(1)}$.
For $t\geq 0$ we obtain $\tilde \xi ^{(2)}_t=\xi ^{(1)}_t$
by means of the Markovian property for
$W_{-\infty }$. For $t\leq 0$
the representation
$W_{-\infty }f=W_{-s}f,\ f\in H^{\xi ^{(1)}}_{[-s},\ s\geq 0$,
gives us $\tilde \xi ^{(2)}_t=W_t\xi ^{(1)}_t=\xi ^{(2)}_t,\ t\leq 0$.

$\Box$

\section {Processes with noncorrelated increments.
The Wold decomposition}

Given a stationary stochastic process $\xi $, there exists the
group of unitaries $U=(U_t)_{t\in {\mathbb R}}$
in the Hilbert space with the inner product defined by
the formula $(\cdot ,\cdot )={\mathbb E}(\cdot \overline \cdot )$
such that $\xi _t=U_t\xi _0,\ t\in {\mathbb R}$.
Remember that the process $\xi $ is said to be nondeterministic
if $\wedge _{t\in \mathbb R}H^{\xi }_{t]}=0$ and deterministic if
$\wedge _{t\in \mathbb R}H^{\xi }_{t]}=H^{\xi }$. It is
evidently that there exist processes which are not nondeterministic
nor deterministic as well.
Given a stationary process $\xi $, it is uniquely defined
the decomposition
$\xi =\xi ^{(1)}\oplus \xi ^{(2)}$, where
$\xi ^{(1)}$ and $\xi ^{(2)}$ are deterministic and nondeterministic
processes correspondingly such that
$\xi ^{(1)}$ and $\xi ^{(2)}$ have noncorrelated ncrements.
In its turn, given a nondetermnistic process
$\xi ^{(2)}$, it is uniquely defined the Wold decomposition
$$
\xi _t^{(2)}=\int \limits _{-\infty }^{t}c(t-s)\zeta (ds),
\eqno (W)
$$
where $\zeta (ds)$ is a noncorrelated measure, such that
${\mathbb E}|\zeta (ds)|^2=ds,$ ${\mathbb E}(\zeta (\Delta )
\overline \zeta (\Delta '))=0$ for all measurable disjoint
sets $\Delta $ and $\Delta '$ (see \cite {Roz}). Thus,
every stationary process $\xi _t=U_t\xi _0$ uniquely defines
the process with noncorrelated stationary increments
$\zeta _t$, which is a invariant curve with respect to
the group $U$. This process can be named by
"a nondeterministic part" of $\xi _t$.
In the following we shall call the Wold decomposition
the possibility to associate with the fixed stationary
process the process with stationary noncorrelated increments
which is its nondeterministic part in the sense given above.
Let $V=(V_t)_{t\in \mathbb R_+}$ be a strong continuous semigroup
of nonunitary isometrical operators in a Hilbert space $H$.
In functional analysis the Wold decomposition is a decomposition
of the form $H=H^{(1)}\oplus H^{(2)}$, where subspaces
$H^{(1)}$ and $H^{(2)}$ reduce the semigroup $V$ to
a semigroup of unitary operators and a semigroup of
completely nonunitary isometrical operators correspondingly. A
completely nonunitary isometrical operator is characterized by the
property that there is no subspace reducing it to a unitary
operator.  Every strong continuous semigroup consisting of
completely nonunitary isometrical operators is unitary equivalent to
its model, which is the semigroup of right shifts $S=(S_t)_{t\mathbb
R_+}$ in the Hilbert space $L^2(\mathbb R_+,{\cal K})$ defined by the
formula $(S_tf)(x)=f(x-t),\ x>t,\ (S_tf)(x)=0,\ 0\leq x\leq t$.
Remember that a deficiency index of the generator
$d=s-\lim \limits _{t\to 0}\frac {V_t-I}{t}$ of
the strong continuous semigroup $V$ is a number
of linear independent solutions to
the equation $d^*f=-f$. The Hilbert space of
values $\cal K$ has the dimension equal to
the deficiency index of the generator of
$V$ (see \cite {Nik}). In the following we shall call
the defciency index of the generator an index of the semigroup.
Notice that every semigroup of completely nonunitary
isometrical operators $V$ with the index
$n>0$ determines $n$ noncorrelated processes
$\xi ^{(i)},\ 1\leq i\leq n,$ with noncorrelated increments
such that $\xi ^{(i)}_{t+s}=\xi ^{(i)}_t+V_t\xi ^{(i)}_s,\
s,t\geq 0$.
In the model case of the semigroup of right shifts
$S$ in $L^2(\mathbb R_+,{\cal K})$ the processes of such type
can be constructed in the following way.
Choose the orthonormal basis of the space
$\cal K$ consisting of the elements $e_i,\ 1\leq i\leq n\leq +\infty
,$ and put $\xi ^{(i)}=e_i\otimes \chi _{[0,t]}$, where $\chi
_{[0,t]}$ is an indicator function of the interval
$[0,t]$. We shall investigate a behaviour of the processes
with noncorrelated increments with respect to "perturbations"
by the Markovian cocycles we introduced in the previous part.

{\bf Proposition 2.1.} {\it Let $\xi $ and $W$ be the process
with noncorrelated increments and the Markovian cocycle
correspondingly.
Then $\xi _t'=W_t\xi _t,\ t\leq 0,\
\xi _t'=\xi _t,\ t>0,$ is a process with noncorrelated
increments.}

Proof.

It follows from Proposition 1.3 that the Markovian isometrical
operator
$W_{-\infty }=s-\lim \limits _{t\to +\infty }W_{-t}$
connects the perturbed process with the initial process
by the formula
$\xi _t'=W_{-\infty }\xi _t,\ t\in {\mathbb R}$.
$\Box $

Let the process with noncorrelated increments
$\xi =(\xi )_{t\in \mathbb R}$ is invariant with
respect to the group of unitaries $U=(U_t)_{t\in \mathbb R}$.

{\bf Proposition 2.2.} {\it The restriction
$V_t=U_{-t}|_{H^{\xi}_{0]}},\ t\geq 0,$ determines
a semigroup of completely nonunitary isometrical operators
with the unit index
in the Hilbert space $H_{0]}^{\xi }$.}

Proof.

Every semigroup of completely nonunitary isometrical operators
with the unit index is unitarily equivalent to its model which
is the semigroup of right shifts $S=(S_t)_{t\geq 0}$
acting in the Hilbert space $L^2({\mathbb R}_+)$
by the formula $(S_tf)(x)=f(x-t),\ x>t,\
(S_tf)(x)=0,\ 0\leq x\leq t$ (see \cite {Nik}).
Define a continuous curve $\eta =(\eta _t)_{t\geq 0}$ in
$L^2({\mathbb R}_+)$ such that
$\eta _t(x)=1,\ 0\leq x\leq t,\ \eta _t(x)=0,\ x>t$.
Linear combinatons of the elements of the curve
$\eta $ form a dense set in the space
$L^2({\mathbb R}_+)$ and
$\eta _{t+s}=\eta _t+S_t\eta _s,\ s,t\geq 0$.
Notice that the stationarity and the orthogonality
of the curve increments imply
$||\xi _t||^2=||\sum \limits _{i=1}^n(\xi _{i\frac
{t}{n}}-\xi _{(i-1)\frac {t}{n}})||^2=n||\xi _{\frac {t}{n}}||^2$ and
$\frac {||\xi _t||^2}{t}=\sigma ^2=const$.  Thus
the curve $\xi $ can be represented in the form $\xi _t=
\sigma ^2\mu _t$, where the measure $d\mu _t$ with values
in $H$ satisfies the conditions
$\mu (\Delta _1)\perp \mu (\Delta _2)$ for
disjoint measurable sets
$\Delta _1,\Delta _2\subset
{\mathbb R}$ ¨ $||d\mu _t||^2=dt$.
Define a unitary operator
$W:H^{\xi }_{0]}\to L^2({\mathbb R}_+)$ by the formula $W\mu_t=\eta
_t$, then $V_t=W^*S_tW,\ t\geq 0$.  $\Box $

Let $V=(V_t)_{t\geq 0}$ be the semigroup of completely
nonunitary isometrical operators in
$H^{\xi }_{0]}$ defined in Proposition 2.2.

{\bf Proposition 2.3.} {\it The Markovian cocycle $W$ determines
a semigroup of isometrical operators $\tilde V=(W_{-t}V_t)_{t\in
\mathbb R_+}$ in $H^{\xi }_{0]}$ with the unit index.
The Wold decomposition $H^{\xi }_{0]}=H^{(1)}\oplus H^{(2)}$
associated with the semigroup $\tilde V$
can be done by the condition
$H^{(2)}=W_{-\infty }H^{\xi }_{0]}$, where $W_{-\infty }=s-\lim
\limits _{t\to +\infty }W_{-t}$.}

Proof.

It follows from Propositions 1.2 and 1.3 that there exists
the limit $W_{-\infty }=
s-\lim \limits _{t\to +\infty }W_{-t}$ such that
the Markovian isometrical operator $W_{-\infty }$ satisfies
the condition
$W_t\xi _t=W_{-\infty }\xi _t,\ t\leq 0$.
Let $H^{(1)}$ be a subspace of $H^{\xi }$
defined by the condition of orthogonality to all elements
$W_{-\infty }\xi _t,\ t\in {\mathbb R}$. This subspace is
invariant with respect to the action of the group of
unitaries $W_tU_t,\ t\in {\mathbb R}$.
On the other hand, by means of the Markovian property for
$W_{-\infty }$, the subspace
$H^{(1)}$ is orthogonal to all elements
$W_{-\infty }\xi _t=\xi _t,\ t\geq 0$, and, therefore,
$H^{(1)}\subset H^{\xi }_{0]}$. Thus, we have proved that
the subspace $H^{(1)}\subset H^{\xi }_{0]}$ is invariant
with respect to the group of unitaries
$(W_tU_t)_{t\in {\mathbb R}}$. Hence, the restriction
$V_t|_{H^{(1)}}=W_{-t}U_{-t}|_{H^{(1)}},\ t\geq 0$, consists of
unitary operators.
The subspace $H^{(2)}\subset H^{\xi }$ is determined by
the condition of orthogonality to
$H^{(1)}$, which is $H^{(2)}=
W_{-\infty }H^{\xi }$ by the definition of $H^{(1)}$.
In this way, the restriction $V|_{H^{(2)}}$ is
a semigroup of completely nonunitary isometrical operators
by means of Proposition 2.2. $\Box $

\section {A model of the Markovian cocycle for the process
with noncorrelated increments}

Let $S=(S_t)_{t\in \mathbb R}$ be a flow of shifts on the line
acting by the formula $(S_t\eta)(x)=\eta(x+t),\ x,t\in
{\mathbb R},\ \eta\in H=L^2({\mathbb R})$.
The group $S$ is naturaly associated with
the stochastic process $\xi _t=\chi _{[-t,0]}$
with noncorrelated increments such that
$\xi _{t+s}=\xi _t+S_t\xi _s,\ s,t\in \mathbb R$.
Let the subspace $H_{t]}$ consist of functions $f$ with
the support $suppf\subset [-t,+\infty )$. Then
$H_{t]}$ is generated by the increments $\xi _s-\xi _r,\ s,r\leq t$.
The restrictions $T_t=S_{-t}|_{H_{0]}},\ t\geq 0,$ form
the semigroup of right shifts $T$. Every
invariant subspace $T_t{\cal V}\subset {\cal V},\ t\geq 0,$
can be described as an image of the isometrical operator $M_{\Theta
},\ {\cal V}=M_{\Theta } H_{0]}$, where $M_{\Theta }={\cal
L}^{-1}\Theta {\cal L}$. Here $\cal L$ is the Laplace
transformation and
$\Theta $ is a multiplication operator by
$\Theta $ which is an inner function in the semiplane
$Re\lambda \geq 0$. Remember that a function $\Theta
(\lambda )$ is said to be inner if it is analitical in
the semiplane $Re\lambda \geq 0$ and its module equals one
on the imaginary axis. (see \cite {Nik}).  Denote
$P_{[0,t]},P_{[t,+\infty )},P_{{\cal V}}, P_{{\cal V}^{\perp}}$
orthogonal projections on the subspaces of functions with
the support belonging to the segment
$[0,t]$ and the semi-interval $[t,+\infty )$ and
on ${\cal V}$ and ${\cal V}^{\perp}$ correspondingly.

Proposition 2.3 shows that arbitrary perturbation
of the group of shifts associated with the process
with noncorrelated increments by a Markovian cocycle $W$
is completely described by a unitary part
$R=\tilde V|_{H^{(1)}}$ of the semigroup of isometries
$\tilde V=(W_{-t}V_t)_{t\in {\mathbb R}_+}$.
In the following theorem we construct the Markovian cocycles
resulting in the semigroup $R$ which is unitarily equivalent
to one we set. In this way, we introduce a model
describing all Markovian cocycles up to unitary equivalence
of perturbations.

{\bf Theorem 3.1.} {\it Let $R=(R_t)_{t\in \mathbb R_+}$
be a strong continuous semigroup of unitaries in the space
${\cal V}^{\perp}$, where $\cal V$ is invariant
with respect to the semigroup of right shifts $T$.
Then the family of unitary operators $(W_t)_{t\in \mathbb R}$,
defined for negative $t$ by the formula
$$
W_{-t}\eta =
(R_tP_{{\cal V}^{\perp}}S_t-P_{\cal V})P_{[t,+\infty )}\eta+
M_{\Theta }P_{[0,t]}\eta ,\ \eta \in H_{0]},
$$
$$
W_{-t}\eta =\eta ,\ \eta \in H_{[0},\
t\geq 0,
$$
and extended for positive $t$ by the formula
$W_t=S_tW_{-t}^*S_{-t},$ $t\geq 0,$
is a Markovian cocycle such that
$\lim \limits _{t\to +\infty }W_{-t}\eta =M_{\Theta }\eta ,$
$\eta \in H_{0]}$. The semigroup $R$ determines a unitary part
of the semigroup of isometries
$(W_{-t}T_t)_{t\in \mathbb R_+}$ in the space
$H_{0]}$ associated with it accordingly to the Wold decomposition.
}

Proof.

Notice that the projection
$P_{\cal V}$ and
the isometrical operator $M_{\Theta }$
are commuting with the right shifts
$T_t=S_{-t},\ t\geq 0$,
and the projections $P_{[t,+\infty )}$.
Hence, the family of unitary operators
$W_{-t}S_{-t}=
(R_tP_{{\cal V}^{\perp}}-S_{-t}P_{\cal V})
P_{H_{0]}}\eta+ (M_{\Theta }P_{[0,t]}+
P_{H_{[0}})S_{-t},\ t\geq 0,$ forms a semigroup
in $H$.
The restriction $V_t=W_{-t}S_{-t}|_{H_{0]}},\ t\geq 0$,
is a semigroup of nonunitary isometrical operators
in $H_{0]}$ with the Wold decomposition
$H={\cal V}^{\perp }\oplus {\cal V}$. Really,
the restriction $V_t|_{{\cal V}}$ is intertwined with
the semigroup of right shifts by
the isometrical operator $M_{\Theta }$
implementing a unitary map of $H_{0]}$ to
${\cal V}$ such that $V_t|_{\cal V}M_{\Theta }=
M_{\Theta }T_t,\ t\geq 0$. Thus,
$V|_{\cal V}$ and $T$ are unitarily equivalent.
The restriction $V|_{{\cal V}^{\perp}}=R$. Therefore,
$V=(W_{-t}T_t)_{t\geq 0}$ is the semigroup of nonunitary
isometrical operators with the unit index and the
unitary part $R$.
$\Box $

Below, using the model of Markovian cocycle introduced
in Theorem 3.1, we construct the Markovian cocycle
satisfying the property $W_t-I\in
s_2,\ t\in \mathbb R$.
Cocycles of such type can be named inner for further applications
in quantum probability (see \cite {Amo001, Amo01, Amo012}).
In fact, this condition appears, particularily, as the condition
of innerness for the quasifree automorphism of the Fermion algebra,
which is generated by the unitary operator $W$
(see \cite {Mur, Amo001, Amo012}).
It is possible to explain why this condition appears in
the following way. Attempts to define a measure in
a Hilbert space $H$ result in constructing the
measure of white noise on the space $E^*$ involved
in the triple $E\subset H\subset
E^*$, where we denote $E^*$ the adjoint space of linear
functionals on the space $E$ which is dense
in $H$ (see \cite {Hid}).
This situation is realized, particularily, if
$\xi \in H$ are included in the parameter set of
the generalized stochastic process.
Given $\xi \in E,\ x\in E^*$, denote $<x,\xi >$
the corresponding dual product.
Notice that if $x\in H$, then $<x,\xi >=(x,\xi )$
coincides with the inner product in $H$.
Fix a positive bounded operator $R$ in the space $H$.
Then there exists a space
$(\Omega ,\mu )$ with the Gaussian measure $\mu $ such that
$\int \limits _{E^*}e^{i<x,\xi
>}d\mu (x)= e^{-(\xi ,R\xi )},\ \xi \in E\subset H$.
Suppose that a unitary operator $W$ in $H$ maps not
elements of $H$ but the measure $\mu $ to certain
other Gaussian measure $\tilde \mu $ such that
$\int \limits _{E^*
} e^{i<x,\xi>}d\tilde \mu (x)=e^{-(\xi ,W^*RW\xi )}, \ \xi \in
E\subset H$.
It is naturally to ask: when are the Gaussian measures
determined by the operators
$R$ and $W^*RW$ equivalent? The Feldman criterion
(see  \cite {Fel, Gui0}) gives us the condition $(\xi ,R\xi )-(\xi
,W^*RW\xi )=(\xi ,\Delta \xi )$, where $\Delta $ is a hermitian
operator of the Hilbert-Schmidt class.  Thus, $R-W^*RW\in
s_2$ which can be rewritten as $WR-RW\in s_2$. The condition
of the type given above is satisfied for all positive
operators $R$ if
$W-I\in s_2$. Notice that the condition of the Feldman type
appears in \cite
{Ara} as the condition of quasi-equivalence for Gaussian states
on the Boson algebra.

Let the inner function $\Theta$ involved in the condition
of Theorem 3.1 be the Blaschke product
(see \cite {Nik}) constructed from the complex numbers
$\lambda _k,\
1\leq k\leq N\leq +\infty $, in the following way.
By means of \cite {Amo001, Amo002, Amo012} take
real parts such that
$Re\lambda _k<0,\ \sum \limits _{k=1}^N|Re\lambda _k|<
+\infty $, imaginary parts can be chosen arbitrary,
then $\Theta (\lambda )=\prod \limits _{k=1}^N
\frac {\lambda +\overline \lambda _k}{\lambda -\lambda
_k}$.  The Blaschke product $\Theta (\lambda )$ is
a regular analytical function in the semiplane
$Re\lambda >0$ and its module equals one
on the imaginary axis. The corresponding subspace ${\cal V}$
of the Hilbert space $H_{0]}=L^2({\mathbb R}_+)$, which is
invariant with respect to the semigroup of right shifts,
is determined by the condition of orthogonality to all
exponents $e^{\lambda _kx},\
1\leq k\leq N$. Let the functions
$g_k,\ 1\leq k\leq N,$ be obtained by
the successive orthogonalization of the system
$e^{\lambda _kx}$. Then
$(g_k,g_l)=\delta _{kl}$ and $g_k,\ 1\leq k\leq N,$
form a orthonormal
basis of the space ${\cal V}^{\perp}$. Define
a $C_0$-semigroup of unitaries $R=(R_t)_{t\in
\mathbb R}$ by the formula $R_tg_k=e^{iIm\lambda _kt}g_k,\
1\leq k\leq N,\ t\in {\mathbb R}_+$.

{\bf Theorem 3.2.} {\it The Markovian cocycle $W=(W_t)_{t\in {\mathbb
R}}$ associated with the inner function $\Theta $
and with the semigroup of unitaries $R$
constructed above as in Theorem 3.1
is inner, i.e. it satisfies the condition
$W_t-I\in s_2,\ t\in {\mathbb R}$.}

{\bf Corrolary 3.3.} {\it
Perturbing the semigroup of right shifts $T$
by the inner Markovian cocycle, it is possible to
obtain the semigroup of isometrical operators with
the unitary part possesing the pure point spectrum
we put.}

Proof of the corrolary. The semigroup
$R$ is a unitary part of the semigroup
$(W_{-t}T_t)_{t\in {\mathbb R}_+}$ in the space $H_{0]}$
by means of Theorem 3.1. The point spectrum of $R$ consists of
imaginary numbers $iIm\lambda _k,$ $1\leq k\leq N\leq +\infty ,$
which can be chosen arbitrary. $\Box $

Proof of Theorem 3.2.

Notice that $W_{-t}\eta -\eta =0,\ \eta \in {\cal H}_{0]}^{\perp},\
t>0$. Thus, we need to prove a convergence of the series
$\sum \limits _{i=1}^{+\infty }||W_{-t}\eta _i-\eta _i||^2$
for a orthonormal basis $(\eta _i)_{i=1}^{+\infty
}$ of the space ${\cal H}_{0]}$. Represent $W_{-t}-I|_{{\cal H}_{0]}}$
as a sum of two parts such that
$W_{-t}\eta -\eta =
(R_tP_{{\cal V}^{\perp}}S_t-P_{{\cal V}^{\perp}})P_{[t,+\infty )}\eta+
(M_{\Theta }-I)P_{[0,t]}\eta ,\ \eta \in H_{0]},\ t>0$, and
prove a convergence of the serieses associated with these parts,
i.e.
$$
\sum \limits _{k=1}^{+\infty }||(R_tP_{{\cal V}^{\perp}}S_t-
P_{{\cal V}^{\perp}})P_{[t,+\infty )}\eta _k||^2<+\infty
\eqno (r1)
$$
and
$$
\sum \limits _{k=1}^{+\infty }||(M_{\Theta }-I)P_{[0,t]}\eta _k||^2
<+\infty .
\eqno (r2)
$$
To check $(r1)$ it is sufficiently to prove a convergence of the
series
$\sum \limits _{k=1}^{N}||(R_t-P_{{\cal V}^{\perp}}S_t^*
P_{{\cal V}^{\perp}})g_k||^2$, where the functions $(g_k)_{k=1}^{N}$,
forming the orthonormal basis of the space
${\cal V}^{\perp}$ are obtained by a successive orthogonalization
of the exponents $e^{\lambda _kx}$. One can represent $R_t$ as a sum
of $R_t^{(1)}$ and $R_t^{(2)}$, where $R_t^{(1)}g_k=e^{\lambda
_kt}g_k$, $R_t^{(2)}g_k= (e^{iIm\lambda _kt}-e^{\lambda _kt})g_k,\
1\leq k\leq N,\ t>0$.  Then the series $\sum \limits
_{k=1}^{N}||(R_t^{(1)}- P_{{\cal V}^{\perp}}S_t^*P_{{\cal V}^{\perp
}})g_k||^2$ converges by the theorem on the triangulation of
the truncated shift (see \cite {Nik}).  Notice that $\sum \limits
_{k=1}^N|e^{iIm\lambda _kt}-e^{\lambda _kt}|= \sum \limits
_{k=1}^N|1-e^{Re\lambda _kt}|$ and the last series converges
because $\sum \limits _{k=1}^N|Re\lambda _k| <+\infty $.
Therefore, $R_t^{(2)}\in s_1$ (the first Schatten class) such that
$$
\sum \limits _{k=1}^N||(R_t- P_{{\cal V}^{\perp}}S_t^*P_{{\cal
V}^{\perp }})g_k||^2\leq \sum \limits _{k=1}^N||(R_t^{(1)}- P_{{\cal
V}^{\perp}}S_t^*P_{{\cal V}^{\perp }})g_k||^2+
$$
$$
\sum \limits
_{k=1}^N||R_t^{(2)}g_k||^2<+\infty .
$$
To prove $(r2)$ it is sufficiently to find
the set of functions $f_k(x),$ $k\in {\mathbb Z},$ which is
dense in the space $H_{[0,t]}=P_{[0,t]}H$ with the property
that there exists the bounded operator $V$ with the bounded inverse
$V^{-1}$ such that the set of functions $(Vf_k)_{k=-\infty }^{+\infty
}$ forms an orthonormal basis in $H_{[0,t]}$ and the following series
converges, $\sum \limits _{k=-\infty }^{+\infty } ||(M_{\Theta
}-I)f_k||^2<+\infty $. The set of functions
$(f_k)_{k=-\infty}^{+\infty }$ satisfying the property given above is
called a Riesz basis of the space $H_{[0,t]}$. A canonical example of
the Riesz basis is the set of exponents $f_k=e^{\mu _kx},\ x\in
[0,t],$ with the indicators $\mu _k=-\frac {1}{2|k|}+ i\frac {2\pi
k}{t},\ k\in {\mathbb Z}$.  Put $f_k(x)=0$ for $x\notin [0,t]$. Then
$f_k=f_k^{(1)}-f_k^{(2)}$, where $f_k^{(1)}=e^{\mu _kx},\
x\in {\mathbb R}_+$, $f_k^{(2)}=e^{\mu _kx},\ x\geq t$,
and $f_k^{(1)}(x)=0,\ x\in {\mathbb R}_-,\ f_k^{(2)}(x)=0,\
x<t$. Using the Parseval equality for the Laplace transformation
and taking into account that the Blaschke product $\Theta $ is
an isometrical operator, we obtain
$$
\sum \limits _{k\in {\mathbb Z}}||(M_{\Theta }-I)f_k||^2\leq
2\sum \limits _{k\in {\mathbb Z}}(||(M_{\Theta }-I)f_k^{(1)}||^2+
||(M_{\Theta }-I)f_k^{(2)}||^2)=
$$
$$
\frac {1}{\pi }\sum \limits _{k\in {\mathbb Z}}(||(\Theta -I)\tilde
f_k^{(1)}||^2+||(\Theta -I)\tilde f_k^{(2)}||^2)
=
$$
$$
\frac {2}{\pi }\sum \limits _{k\in {\mathbb Z}}\{
(||\tilde f_k^{(1)}||^2-Re(\Theta \tilde f_k^{(1)},\tilde f_k^{(1)}))+
(||\tilde f_k^{(2)}||^2-Re(\Theta \tilde f_k^{(2)},\tilde
f_k^{(2)}))\},
\eqno (B1)
$$
where
$\tilde f_k^{(1)}(\lambda )=\frac {1}{\lambda
\mu _k},\ f_k^{(2)}=e^{\mu _kt}\frac {e^{-t\lambda }}{\lambda -\mu
_k}$ is the Laplace transformation of the functions
$f_k^{(1)}$ and $f_k^{(2)}$. Apply the anlitical functions
techniques, we get
$$
(\Theta \tilde f_k^{(i)},\tilde f_k^{(i)})=\Theta
(-\overline {\mu _k})||f_k^{(i)}||^2,\ i=1,2.
\eqno (B2)
$$
Notice that $ln\Theta (\lambda )=\sum \limits _{k=1}^N(ln(1+\frac {\overline
{\lambda }_k}{\lambda })- ln(1-\frac {\lambda _k}{\lambda }))=
\frac {\sum \limits _{k=1}^NRe\lambda _k}{\lambda }+o(\frac
{1}{\lambda })$.
Hence,
$$
\Theta (\lambda )=1+\frac {\sum \limits
_{k=1}^{+\infty }Re\lambda _k}{\lambda }+o(\frac {1}{\lambda }).
\eqno (B3)
$$
Substituting in $(B1)$ the formulas $(B2)$ and $(B3)$, we obtain
$$
\sum \limits _{k\in {\mathbb Z}}||(M_{\Theta }-I)f_k||^2\leq
C_1\sum \limits _{k\in {\mathbb Z}}(1-Re\Theta (-\overline {\mu
_k}))\leq
$$
$$
C_2\sum \limits _{k\in {\mathbb Z}}\frac {\frac
{1}{2|k|}}{\frac {1}{4k^2}+ \frac {4\pi ^2k^2}{t^2}}\leq
C_3\frac {1}{|k|^3}<+\infty ,
$$
where $C_1,C_2$ and $C_3$ are some positive constatnts.
$\Box $

\section {Processes with independent increments in
classical and quantum probability. Kolmogorov
flows}

Denote ${\cal L}_s(H)$ and $\sigma _1(H)$ the sets of linear
hermitian and positive unit-trace operators in a Hilbert
space $H$.
In the quantum probability theory the elements $x\in {\cal L}_s(H)$
and $\rho \in \sigma _1(H)$
are called random variables or observables and
states of the system correspondingly.
Consider the spectral decomposition
$x=\int \lambda dE_{\lambda }$ of the random variable $x\in {\cal
L}_s(H)$, where $E_{\lambda }$ is a resolution of the identity in
$H$. Then a probability distribution of $x$ in the state $\rho
\in \sigma _1(H)$ is defined by the formula
$P(x<\lambda )=Tr\rho E_{\lambda }$.
Thus, the expectation of $x$ in the state $\rho
$ can be calculated as ${\mathbb E}(x)=Tr\rho x$ (see \cite {Hol0,
Hol1}). Notice that the classical random variables from
$L^{\infty}(\Omega )$ can be considered as linear operators
multiplying by the function in the Hilbert space
$H=L^2(\Omega )$, where $\Omega $ is some probability space.
A quantum stochastic process (in the narrow sense of the word)
is a strong continuous family of operators
$x_t\in {\cal L}_s(H),\ t\in
\mathbb R$, "quantum observables". We shall call the quantum
stochastic process by a process with stationary increments
if there exists a ultraweek continuous one-parameter group
of *-automorphisms $\alpha _t,\ t\in \mathbb R,$ of the
algebra of all bounded operators in $H$ such that
$x_{t+s}=x_t+\alpha _t(x_s),\ s,t\in \mathbb R$.
We do not suppose that the operators $x_t$ are bounded
but assume that the action of $\alpha =(\alpha _t)_{t\in \mathbb R}$
is correctly defined on $x=(x_t)_{t\in \mathbb R}$.
The quantum stochastic process with stationary increments
$x_t$ is called a stationary process if
there exists $x\in {\cal L}_s(H)$ such that $x_t=\alpha _t(x)-x,\ t\in
{\mathbb R}$. As in the case of classical stochastic processes,
the quantum stochastic process which is continuous in the square mean
determines a continuous curve in the Hilbert space
with the inner product defined by the expectation.
This curve we shall denote $[x]= ([x_t])_{t\in \mathbb R}$.
As we identified continuous curves in a Hilbert space
with classical stochastic processes, we obtain that
given a quantum stochastic process $x$ can be associated
with the classical stochastic process $\xi =[x]$.
If $x$ is a process with stationary increments or a stationary
process and the state $\rho \in \sigma _1(H)$ is invariant
with respect to the action of the group $\alpha $, i.e. $Tr(\rho
\alpha _t(a))=Tr(\rho a),\ t\in \mathbb R$, where $a$ is
arbitrary linear combination of the random variables
$x_t,\ t\in
\mathbb R$, then the classical stochastic process $\xi =[x]$ is
also stationary. The group of unitaries $U$ shifting
the increments $\xi $ in time is defined by the formula $(U_t\xi
_s,\xi _r)=Tr(\rho \alpha _t(x_s)x_r),\ s,t,r\in \mathbb R$.
Because the notion of expectation in the quantum probability theory
plays a major role, we shall use it to define an independence
in the classical probability theory also for convenience
to pass to quantum probability in the following.
The classical stochastic process $x=(x_t)_{t\in \mathbb
R}$ is said to be a
process with independent increments or a Levy process if
the following identity holds,
$$
{\mathbb E}(\phi _1(x _{t_1}-x _{s_1}) \phi
_2(x _{t_2}-x _{s_2})\dots \phi _n(x _{t_n}-x _{s_n}))= \prod \limits
_{i=1}^n{\mathbb E}(\phi _i(x _{t_i}-x _{s_i}))
\eqno (I1)
$$
for arbitrary choise of functions $\phi _i\in L^{\infty }$
and disjoint intervals $(s_i,t_i)$.  There are different
approaches to a definition of the Levy processes in quantum
probability (see \cite {Hol0, Hol1}). In any way, besides that
the condition $(I1)$ must be satisfied, they are involved certain
additional conditions concerning the algebraic structure of the
process $x$. One of them is a commutativity for increments, that is
$$
[x_{t_1}-x_{s_1},x_{t_2}-x_{s_2}]=0
\eqno (I2)
$$
for disjoint intervals $(s_i,t_i)$.
Notice that for the quantum processes with independent increments
it is possible to define the representation of
the Levy-Hinchin type (see \cite {Hol2, Hol3}).
For convenience we shall bring the definition of the Kolmogorov
flow from the Introduction.
The flow $\{T_t,\ t\in {\mathbb R}\}$ on the probability space
$(\Omega ,{\cal M},\mu )$ is said to be a Kolmogorov flow
(see \cite {Kol3}) if there exists a $\sigma $-algebra
of events ${\cal M}_{0]}\subset
{\cal M}$, such that $T_t{\cal M}_{0]}={\cal M}_{t]},$
$$
{\cal M}_{s]}\subset {\cal
M}_{t]},\ s<t, \eqno (K_1)
$$
$$ \cup _t{\cal M}_{t]}={\cal M},
\eqno (K_2)
$$
$$
\cap _t{\cal M}_{t]}=\{\emptyset ,\Omega \}.
\eqno (K_3)
$$
In the formula $(K_3)$ we mean that the intersection
of $\sigma $-algebras contains only two events which are
the empty set and the all space $\Omega $.
Let ${\cal M}_{t]}$ be generated by the events associated
with the stationary stochastic process $\xi _s,\ s<t$.
The investigation of the conditions which lead to
the Kolmogorov flow generated by $\xi $, it is given in
\cite {Ibr}. Particularily, the flow of the Wiener process
is a Kolmogorov flow (see \cite {Hid}).
Notice that to obtain the Kolmogorov flow from the stochastic
process, it is not nessesary to claim the independence of
the increments. Consider the quantum stochastic process
with stationary increments
$x=(x_t)_{t\in \mathbb R}$. In the quantum probability
theory a role of the $\sigma $-algebras of events is played
by the von Neumann algebras generated by the quantum random
variables. In this way, the conditional expectation is
a completely positive projection on the von Neumann algebra
(see \cite {Hol0, Hol1}).  Let ${\cal M}_{t]}=\{x_s,\ s<t\}''$ and
${\cal M}_{[t}= \{x_s,\ s>t\}$ be the von Neumann algebras
generated by the past before the moment $t$ and the future after
the moment $t$ of the quantum stochastic process with stationary
increments $x$. Here
$A'$ denotes the set of all bounded operators in $H$
which are commuting with the operators (not bounded in general)
from the set $A$.
It follows from the stationarity of the increments of
the process $x$ that
${\cal M}_{t]}=\alpha _t({\cal M}_{0]}),\ t\in {\mathbb R},$
for the group $\alpha $ shifting the increments.
It is naturally to call the group of automorphisms $\alpha $
a Kolmogorov flow on the von Neumann algebra $\cal M$ if
the conditions $(K_1)$, $(K_2)$, $(K_3)$ are satisfyied, where
in the condition $(K_3)$ the trivial $\sigma
$-algebra $\{ \emptyset ,\Omega \}$ is replaced to the trivial
von Neumann algebra $\{{\mathbb C}{\bf 1}\}$
containing only operators which are multiple to the identity.
Thus we obtain the condition of the algebraic Kolmogorov flow
$\alpha $ (see \cite {Emc}).

For classical random variables $\xi $ with ${\mathbb E}\xi =0$,
the condition $D\xi ={\mathbb E}\xi ^2=0$ means that $\xi $
equals zero almost surely.
The situation changes under a transformation to the quantum
case. For a quantum random variable $x$
with ${\mathbb E}x=0$ the equality $Dx={\mathbb E}x^2=
Tr\rho x^2=0$ doesn't imply $x=0$. The property of this type
characterizes the expectation ${\mathbb E}={\mathbb E}_
{\rho }$. If the situation appears, that is
$Tr\rho x=0$ implies $x=0$ for all positive operators
$x$ belonging to certain algebra $\cal M$, then
$\rho $ is said to determine a faithful state
${\mathbb E}(\cdot )=
Tr\rho \cdot $ on $\cal M$.

{\bf Proposition 4.1.}{\it Let the expectation
${\mathbb E}(\cdot )=Tr\rho \cdot $ determine a
faithful state on the von Neumann algebra $\cal M$ generated
by the stationary increments of the quantum stochastic process
satisfying the condition $(I1)$.
Then the group of automorphisms $\alpha $ shifting
the increments in time is the Kolmogorov flow.}

Proof. Let $x\in \cap _t{\cal M}_{t]}$,
then given $y\in \cup _t{\cal M}_{[t}={\cal M}$
the condition $(I1)$ implies ${\mathbb E}((x-{\mathbb Ex})y)=
{\mathbb E}(x-{\mathbb E}(x){\bf 1}){\mathbb E}(y)=
0=Tr\rho (x-{\mathbb E}(x))y$. Put
$y=x-{\mathbb E}(x)$. Because the state $\rho $ is
faithful, it follows that the identity
$Tr\rho (x-{\mathbb E}(x))^2=0$ implies
$x={\mathbb E}(x)\in \{{\mathbb C}{\bf 1}\}$.
$\Box $

In the following part we shall give the important example
of the algebraic Kolmogorov flow, which complements
the example of Proposition 4.1 in some sense.

\section {Cohomology of groups as a language describing
the perturbations of the space of all functionals from
stochastic process}

Below we remember some notion of the cohomology of groups
(see, f.e., \cite {Bra, Gui}).
Define certain action $\alpha $ of the real line
$\mathbb R$ as an additive group on the algebra
$\cal M$. The element $I\in Hom({\mathbb
R}^k, {\cal M})$ is said to be (additive)
$k-\alpha $-cocycle if the following identity
is satisfied,
$$
\alpha _{t_1}I(t_2,t_3,\dots ,t_{k+1})-
I(t_1+t_2,t_3,\dots ,t_{k+1})+\dots
$$
$$
+(-1)^i
I(t_1,\dots ,t_{i-1},t_i+t_{i+1},t_{i+2},\dots ,t_{k+1})+\dots
$$
$$
+(-1)^{k+1}I(t_1,\dots ,t_k)=0,
$$
$\ t_i\in {\mathbb R},\ 1\leq i\leq k+1$;
the $k-\alpha $-cocycle $I$ is said to be a coboundary
if there exists the element
$J\in Hom({\mathbb R}^{k-1},{\cal M})$ such that
$$
I(t_1,\dots ,t_k)=\alpha _{t_1}(J(t_2,\dots ,t_k))-
J(t_1+t_2,\dots ,t_k)+\dots
$$
$$
+(-1)^iJ(t_1,\dots ,t_{i-1},t_i+t_{i+1},t_{i+2},\dots ,t_k)+\dots
+(-1)^{k+1}
J(t_1,\dots ,t_{k-1}),
$$
$t_i\in {\mathbb R},\ 1\leq i\leq k$.
Denote $C^k$ and $B^k$ the sets of all additive
$k-\alpha$-cocycles and
$k-\alpha $-coboundaries correspondingly. Then
$H^k=H^k(\alpha )=C^k/B^k$ is called a $k$-th group of
cohomologies of $\alpha $ with values in $\cal M$.
Define a cohomological multiplication
$\cup:H^k\times H^l\to H^{k+l}$ by the formula
$$
(I_k\cup I_l)(t_1,\dots ,t_{k+l})=I_k(t_1,\dots ,t_k) \alpha _{t_1+\dots
+t_k}(I_l(t_{k+1},\dots ,t_{k+l})),
$$
$I_k\in C^k,\ I_l\in C^l,\ t_i\in
{\mathbb R},\ 1\leq i\leq k+l$.
The group $
H=H(\alpha )=\oplus _{i=1}^{+\infty }H^i$
is a graded ring with respect to the multiplication
$\cup$.
Notice that a group of $0$-cohomologies was omitted
because we do not need it in the following.
Now let $\alpha $ and $\alpha '$ be two actions
on the algebra $\cal M$. We shall say that
the rings $H(\alpha )$ and $H(\alpha ')$ are isomorphic if
there is a one-to-one correspondence
$w:H(\alpha )\to H(\alpha ')$ mapping each $H^k(\alpha )$ to
$H^k(\alpha ')$ and the action of $\alpha $ to the
action of $\alpha '$.

The one-parameter family of automorphisms
$w=(w_t)_{t\in \mathbb R}$ of the algebra $\cal M$
is said to be {\it a (multiplicative) $\alpha $-cocycle}
if the following condition holds,
$w_{t+s}=w_t\circ \alpha _t\circ w_s
\circ \alpha _{-t},\ s,t\in {\mathbb R},\
w_0=Id$.
We shall call a multiplicative $\alpha $-cocycle by {\it Markovian}
with repect to certain set $\cal I$ of additive
$1-\alpha $-cocycles $I(t)\in {\cal I}$ if the condition
$w_t(I(t+s)-I(t))= I(t+s)-I(t),\ t,s\geq 0$, is satisfied.
The following proposition can be considered as some "abstract
generalization" of the properties of the Markovian perturbations
we have described in the previous parts.

{\bf Proposition 5.1.} {\it Let the ring $A\subset H(\alpha )$
be generated by a set of additive $1-\alpha $-cocycles
$I(t)\in {\cal I}$ and the multiplication
$\cup$. Then given a multiplicative $\alpha $-cocycle $w$ which
is Markovian with respect to $\cal I$ determines the homomorphism
of $A$ into $H(\alpha ')$, where the group $\alpha '$ is defined
by the formula $\alpha '_t=w_t\circ \alpha _t,\ t\in \mathbb R$.
The image $A'$ of this homomorphism is a ring generated by
the set of additive $1-\alpha '$-cocycles $I'(t)=w_t(I(t)),\ t\leq
0,\ I'(t)=I(t),\ t>0,\ I(t)\in {\cal I}$.}

Proof. Check that $I'(t)=w_t(I(t))$ satisfies
the condition for a $1-\alpha '$-cocycle.
Fix $s,t>0$ and notice that
$w_{-s}(I(-s))$ $=w_{-s-t+t}(I(-s))=w_{-s-t}\circ \alpha _{-s-t}\circ
w_t\circ \alpha _{s+t}(I(-s))=w_{-s-t}\circ \alpha _{-s-t}\circ
w_t(I(t)-I(s+t))=w_{-s-t}\circ \alpha _{-s-t}(I(t)-I(s+t))=
w_{-s-t}(I(-s))$ by means of the cocycle condition and
the Markovian property for $w$.
Therefore,
$
w_{-t-s}(I(-t-s))=
w_{-t-s}(I(-t)+\alpha _{-t}(I(-s)))=w_{-t-s}(I(-t))+
w_{-t-s}\circ \alpha _{-t}(I(-s))=w_{-t}(I(-t))+
w_{-t}\circ \alpha _{-t}\circ w_{-s}(I(-s))
$
$
=I'(-t)+\alpha '_{-t}(I'(-s)).
$
Using the Markovian property of $w$ we obtain
$I'(t)=I(t)=-\alpha _t(I(-t))=-w_t\circ \alpha _t\circ w_{-t}(I(-t))=
-\alpha '_t(I'(-t)),\ t>0$, because $w_t\circ \alpha _t\circ
w_{-t}\circ \alpha _{-t}=Id$ due to the cocycle condition for $w$.
$\Box $

Below we give the examples showing how the language of the theory
of cohomology of groups can describe the set of all functionals
from the classical and quantum stochastic processes with the
stationary independent increments.

\subsection {The Wiener process}

Consider the Wiener process
$\{B(t),\ t\in {\mathbb R}\}$ implemented on
the probability space $(\Omega ,
\mu )$.
Let $S=(S_t)_{t\in \mathbb R}$ be the group of
transformations shifting the increments of the
process
in time. Then the Wiener process satisfies
the condition of the additive
$1-S$-cocycle, i.e. $B(t+s)=B(t)+S_t(B(s)),\
s,t\in \mathbb R$.
Consider the ring of cohomologies $H(S)$ for the group
$S$ with values in $L^{\infty }(\Omega )$. Then
the Wiener process $B(t)$ generates the subring
$A\subset H(S)$. Denote $\cal H$ the Hilbert space
of all $L^2$-functionals from the Wiener process.
As it is known, one can define for $\cal H$
the Wiener-Ito decomposition ${\cal H}=\oplus _{i=0}^{+\infty }
{\cal H}_i$ in the orthogonal sum of spaces formed by
polinomials of increasing degrees (see, f.e.,
\cite {Hid}).
Notice that the representation
of the graded ring $A$ as the sum of cohomologies of
all degrees
$A=\oplus _{i=1}^{+\infty }
A_i$ is a cohomological analog of the Wiener-Ito
decomposition.
Take a function
$a(x)\in L^2_{loc}({\mathbb R})$
and detrmine a one-parameter family of linear maps $w=(w_t)_{t\in
\mathbb R}$ acting on
$A$ by the formula
$$
w_t(B(t+s)-B(t))=B(t+s)-B(t),\ s\geq 0,
$$
$$
w_t(B(s))=exp \{-\frac {1}{2}\int \limits
_{0}^{s}a(x)dB(x)-\frac {1}{4}\int \limits _0^t|a(x)|^2dx\}\cdot
$$
$$
(B(s)+\int \limits_{0}^ta(x)dB(x)),\ s\leq t.
\eqno (W)
$$

Every $w_t$ defines a unitary transformation
in the space of
$L^2$ - functionals of the white noise (see \cite {Hid})
and satisfyies the property for the Markovian cocycle
by the definition. Thus the following Proposition is hold.

{\bf Proposition 5.2.} {\it The family $w$ is a Markovian
cocycle.
}

\subsection {Quantum noises}

Let $G$ be certain Lie algebra with the involution. Then
(see \cite
{Sch}) there exists a one-parameter family of
*-homomorphisms $j:x\to j_t(x)\in {\cal
L(H)}$ mapping each element $x\in G$ to a strong continuous
one-parameter family of operators $j_t(x),\ t\in \mathbb R,$ in
the symmetric Fock space ${\cal H}={\cal F}(L^2( {\mathbb
R},{\cal K}))$ over the one-particle Hilbert space
$L^2({\mathbb R}, {\cal K})$ consisting of functions on
the real line with values in the Hilbert space
$\cal K$.  Every homomorphism $j_t$ preserves
the commutator such that $[j_t(x),j_t(y)]=[x,y],$ $x,y\in G$, and
satisfies the condition of additive
$1-\alpha $-cocycle with respect to the group of automorphisms
$\alpha $ generated by the group of shifts in the space
$L^2({\mathbb R},{\cal K})$.  Notice that
$j$ can be continued to a *-homomorphism of the universal
enveloping algebra ${\cal U}(G)$ determining the quantum Levy
process (see \cite {Sch}).
In the case if $j_t(x)\in {\cal L}_s({\cal H}),\ t\in \mathbb R,$
under fixed $x\in G$, the one-parameter family $x_t=j_t(x),\ t\in
\mathbb R,$ is the quantum stochastic process in the sense
of the definition given in Part 4. Moreover,
$x_t=j_t(x)$ is the process with independent increments, where
the independence means that the conditions
$(I1)$ and $(I2)$ of Part 4 are satisfied.  Nevertheless, we
shall consider the curves $x_t=j_t(x),\ t\in {\mathbb R},$
consisting of non-hermitian operators as well.
The scope of all curves
$x_t=j_t(x),\ x\in G,$ we shall name {\it a quantum noise}
because the increments of all curves $x_t$ are independent.
Notice that the quantum stochastic processes
$j_t(x)$ can play roles of the Wiener and Poisson processes
in the quantum case (see \cite {Par, Hol0,
Hol1} for references and comments). Consider the ring
of cohomologies $A$ generated by additive $1-\alpha$-cocycles
$j_t(x),\ x\in G$.  Notice that the Fock space can be factorized
such that ${\cal F}(L^2({\mathbb R},{\cal K}))={\cal
F}(L^2 ({\mathbb R}_+,{\cal K}))\otimes {\cal F} (L^2({\mathbb
R}_-,{\cal K}))$. In the following we construct the example of
a multiplicative $\alpha $-
cocycle generating the homomorphism of $A$ which uses
this factorization.
Put $w_{-t}(\cdot )=W_t\otimes I\cdot W_t^*\otimes I,\ t\geq 0,$
where the family of unitary operators $W_t,\ t\geq 0,$
in the Fock space ${\cal F}(L^2({\mathbb R}_+,{\cal K}))$
satisfies the quantum stochastic differential equation
constructed in Example 25.17, \cite {Par}, P. 198,
i.e.
$$
dW=W(dA_m^++d\Lambda _{U^{-1}}-dA_{U^{-1}m}-\frac{1}{2}<<m,m>>),
$$
where $A_m^+,A_m$ and $\Lambda _{U^{-1}}$ are the basic processes
of creation, annihilation and number of particles generated by
a function $m\in {\cal K}$ and a unitary operator
$U$ (see \cite {Par}).

Let ${\cal M}_{t]}$ and ${\cal M}_{[t}$ be the von Neumann
algebras generated by the increments of the quantum noise
$j_t(x),\ x\in G,$
before the moment of time $t$ and after the moment of time $t$
correspondingly.  Then for the von Neumann algebra $\cal M$
generated by all operators $j_t(x),\ x\in G,\ t\in {\mathbb R},$
we obtain ${\cal M}=\vee _t{\cal M}_{t]}=\vee _t{\cal M}_{[t}$,
where ${\cal M}_{t]}=\alpha _t({\cal M}_{0]}),\ t\in {\mathbb R}$.

{\bf Proposition 5.3.} {\it Let the von Neumann algebra ${\cal M}$
generated by increments of the quantum noise $j_t(x),\ x \in G,$
be a factor. Then the group of automorphisms $\alpha $
is a Kolmogorov flow on ${\cal M}$.}

Proof.

Take $x\in \cap _t{\cal M}_{t]}\subset {\cal M}$.
The increments of the quantum noise are independent
quantum random variables by the condition of the Proposition.
Hence, the algebra ${\cal M}_{[t}$ belongs to the commutant
of the algebra ${\cal M}_{t]}$. Therefore, the operator
$x$ is commuting with all operators
$y\in \cup _t{\cal M}_{[t}={\cal M}$.
As the von Neumann algebra ${\cal M}$ is a factor, i.e.
${\cal M}\cap {\cal
M}'=\{{\mathbb C}{\bf 1}\}$, we get $x=const{\bf 1}$. $\Box $

Notice that in the applications the algebra $\cal M$
is often the algebra of all bounded operators, i.e.
the factor of type $\rm I$. In this case the condition of
Proposition 5.3 is satisfied. Take into account that the expectation
in the quantum case is often determined by a pure state
$\rho $ such that
${\mathbb E}(x)=(\Omega ,x\Omega )$. Here $\Omega $ is
some vector in the Hilbert space $\cal H$, where
the operators $x$ act. The pure state on the algebra of
all bounded operators in $\cal H$ can not be faithful
(see \cite {Brat}).
Hence, Proposition 5.3 complements Proposition 4.1.

\section {Markovian perturbations of stationary quantum
stochastic processes
}

At first, we shall generalize the defnition of the quantum stochastic
process given in Part 4.  Following to \cite {Acc2} by {\it a quantum
stochastic process} we shall call the one-parameter family of
*-homomorphisms $j=(j_t)_{t\in \mathbb R}$ from certain algebra with
the involution $\cal A$ to the algebra of linear (unbounded in
general) operators $\cal L(H)$ in some Hilbert space $\cal H$. We
shall suppose that $j_0(x)=0,\ x\in {\cal A}$.  For $\cal A$ one can,
in particular, take certain Lie algebra.
Consider the minimal von Neumann algebra
$\cal M$ generated by all operators $j_t(x),\
x\in {\cal A},$ such that ${\cal M}=\{j_t(x),\ x\in {\cal A}\}''$,
where $'$ denotes the commutant in the algebra
$\cal B(H)$ consisting of all bounded operators in
$\cal H$.  For fixed $s,t\in
\mathbb R$ we denote ${\cal M}_{t]},\ {\cal M}_{[t}$ ¨ ${\cal
M}_{[s,t]}$ the von Neumann algebras generated by the operators
$\{j_t(x)-j_r(x),\ r\leq t\},\ \{j_r(x)-j_t(x),\ r\geq t\}$ and
$\{j_r(x)-j_t(x),\ s\leq r\leq t\}$ correspondingly.
Futher, we shall assume that the state
$\rho >0,\ Tr\rho =1,$ determining
the expectation ${\mathbb E}(y)=Tr\rho y$ for
quantum random variables $y\in \cal M$ is fixed.
Moreover, let $\cal M$ be in "the standard form"
with respect to the state $\rho $, that is
the Hilbert space $\cal H$, where $\cal M$ acts in,
is defined by the map $x\to [x]$ from $\cal M$
to the dense set in $\cal H$ such that the inner product
is given by the formula $([x],[y])={\mathbb E}(xy^*),\ x,y\in {\cal
M}$. The quantum stochastic process $j$ is said to be {\it
stationary}
if there exists a group of automorphisms
$\alpha _t\in Aut({\cal M}),\ t\in \mathbb R,$
whose actions are correctly defined on the operators
$j_s(x),\ \cal A$, such that $\alpha
_t(j_s(x))=j_{s+t}(x),\ x\in \cal A$, and the expectation
${\mathbb E}(\cdot)$ is invariant with respect to the group
$\alpha =(\alpha _t)_{t\in {\mathbb R}}$, i.e. ${\mathbb
E}(\alpha _t(y))={\mathbb E}(y),\ y\in {\cal M}$.
The stationary quantum stochastic process is a particular case
of the quantum stochastic process with the stationary increments
generated by *-homomorphisms $j_t$ which determine the curves
$j_t(x)$ being additive $1-\alpha $-cocycles for every fixed
$x\in {\cal A}$. For the example of the quantum stochastic process
with stationary increments can be choosen the quantum noise.
We shall suppose that $\alpha $ has ultra-week continuous orbits,
i.e. all functions
$\eta (\alpha _t(y))$ are continuous in $t$ for
arbitrary $\eta \in {\cal
M}_*,\ y\in {\cal M}$. The one-parameter family of *-automorphisms
$w_t,\ t\in \mathbb
R,$ of the algebra $\cal M$ is said to be {\it a (multiplicative)
Markovian $\alpha $-cocycle} if the following two conditions are
satisfied:

$(i)\ w_{s+t}=w_s\circ \alpha _s\circ w_t\circ \alpha _{-s},\
s,t\in \mathbb R,$

$(ii)\ w_t(x)=x,\ x\in {\cal M}_{[t},\ t\geq 0$.

Notice that if the von Neumann algebra ${\cal M}_{[t}$
be generated by the increments of the quantum noise
$\{j_t(x),\ x\in {\cal A}\}$, the Markovian cocycle $w$
in the sense of the definition given here is a cocycle
being Markovian with respect to the set $\{j_t(x),\
x\in {\cal A}\}$ in the sense of the definition of
Part 5.

{\bf Proposition 6.1.} {\it Let $j_t(x),\ x\in {\cal A},$
and $w$ be the quantum noise and the Markovian cocycle
correspondingly.
Suppose that $w$ doesn't change values of the expectation
${\mathbb E}$ determining the probability distribution of
$j$. Then the one-parameter family of homomorphisms
$\tilde j_t(x)=w_t\circ j_t(x),\ t\leq 0,\ \tilde j_t(x)=j_t(x),\
t>0,\ x\in {\cal A},$ is the quantum noise isomorphic
to the initial one. In particular, the processes
$\tilde j_t(x),\ x\in {\cal A},$ have independent increments.}

Proof.

The quantum noise $j$ is a quantum stochastic process with
independent increments.
Denote $\alpha $ the corresponding group of automorphisms
shifting the increments in time.
Then, due to Proposition 5.1,
$\tilde j_t(x)$ is an additive
$1-\tilde \alpha $-cocycle for each fixed $x\in {\cal A}$,
where the group $\tilde \alpha $ consists of automorphisms $\tilde
\alpha _t=w_t\circ \alpha _t,\ t\in {\mathbb R}$.
By means of the property $(I2)$ guaranteeing
the independence for increments of the quantum noise $j$,
we obtain for the commutator,
$$
[j_{t_1}(x_1)-j_{s_1}(x_1),j_{t_2}(x_2)-j_{s_2}(x_2)]=0
\eqno (K)
$$
for all $x_1,x_2\in {\cal A}$ and disjoint
intervals $(s_1,t_1)$ and $(s_2,t_2)$. Below we shall
prove that this property holds for the process
$\tilde j$ also.
For $s_1,t_1,s_2,t_2$ it takes place because
$\tilde j_r(x)=j_r(x),\ x\in {\cal A},\ r\geq 0$.
The identity for the case where at least
one of $s_1,t_1,s_2,t_2$ is less than zero
can be obtained by applying the automorphisms
$\tilde \alpha _{-r}=w_{-r}\circ \alpha _{-r},\ r>0,$
to the formula $(K)$. In fact, the automorphism
$\tilde \alpha _{-r}$ shifts the increments for
$r$ units backward and the formula
$(K)$ which is true for the process $\tilde j$ with
positive values $s_1,t_1,s_2,t_2$  automatically
appears to be true
for negative values.
Moreover, arguing in the way given above, it
is easy to obtain that all algebraic properties
satisfied for the operators
$j_t(x),\ x\in {\cal A},\
t\in {\mathbb R},$ are also satisfied for
the operators $\tilde j_t(x),\ x\in {\cal A},\
t\in {\mathbb R}$. It remains to check
the condition of independence $(I1)$ of Part 4
for increments of the process $\tilde j$.
But it holds because it is true for the process $j$
and the expectation is invariant with respect to
the action of the cocycle $w$ by the condition.
Thus, the processes $j$ and $\tilde j$ are
isomorphic. $\Box $

Consider the restriction $\beta _t=\alpha _{-t},\ t\in {\mathbb
R}_+,$ to the subalgebra ${\cal M}_{0]}$ which is invariant with
respect to the action of $\alpha _{-t},\ t\in {\mathbb R}_+$.
The unital semigroup $\beta =(\beta _t)_{t\in {\mathbb R}_+}$
consists of the endomorphisms of ${\cal M}_{0]}$ possessing
the property
$\beta ({\cal M}_{0]})\neq {\cal M}_{0]}$ and
has the orbits continuous in the sense that
$\eta (\beta _t(x))$ is a continuous function for
all $\eta \in {\cal M}_{0]*},\
x\in {\cal M}_{0]}$. Hence, $\beta $ is a $E_0$-semigroup
by the definition introduced in \cite {Pow}. If $\alpha $
is a Kolmogorov flow, the $E_0$-semigroup $\beta $ is
a semiflow of Powers shifts, i.e. each
$\beta _t,\ t>0,$ is a Powers shift \cite {Pow} such that
$\cap _{n\in {\mathbb N}}\beta _{tn}({\cal M}_{0]})=\{{\mathbb C}{\bf
1}\},\ t>0$ (see \cite {Bul1}).
Notice that the Markovian cocycle $w$ generates a new group
of automorphisms $\tilde \alpha _t=w_t\circ \alpha _t,\
t\in \mathbb R$, on the von Neumann algebra
$\cal M$ and a new $E_0$-
semigroup $\tilde \beta _t=w_{-t}\circ \beta _t,\ t\in
{\mathbb R}_+,$ on the von Neumann algebra ${\cal M}_{0]}$.
Remember that the conditional expectation of the unital
algebra ${\cal M}$ onto the unital algebra
${\cal N}\subset {\cal M}$ is a completely positive
projection of ${\cal M}$ onto
${\cal N}$ (see, f.e., \cite {Hol0, Hol1}).
Using the techniques \cite {Bul} one can exclude
the maximum subalgebra of the algebra
${\cal M}_{0]}$ such that the restriction of
the semigroup $\tilde \beta $ to it
is a semigroup of automorphisms and there exists
a conditional expectation onto this algebra.
The quantum random variables
$\tilde j_t(x)=w_t\circ j_t(x),\ x\in {\cal A},$ generate
the von Neumann subalgebras
$\tilde {\cal M}=\{\tilde j_t(x),\ x\in {\cal
A},\ t\in {\mathbb R}\}''$ and $\tilde {\cal M}_+=\{\tilde j_{-t}(x),\
x\in {\cal A},\ t\in {\mathbb R}_+\}''$ of the algebras
$\cal M$ and
${\cal M}_{0]}$ correspondingly. Suppose that $\cal M$
is a factor. Then, due to Proposition 5.3, $\alpha $
is a Kolmogorov flow on ${\cal M}$ and by means of Proposition
6.1 the restriction $\tilde \alpha |_{\tilde {\cal M}}$ is also
a Kolmogorov flow. Thus, $\tilde \beta |_{\tilde {\cal
M}_+}$ is a semiflow of Powers shifts.
The Wold decomposition $(W)$ given in
Part 2 allows to uniquely determine the stochastic
process with noncorrelated increments associated with
the stationary process.
In the following theorem we establish the possibility to exclude
a restriction of the group of automorphisms
obtained by a perturbation of the Kolmogorov flow
generated by the quantum noise, which is isomorphic to
the initial Kolmogorov flow.
In this way, our conjecture can be considered as some analog
of the Wold decomposition for the quantum case.

{\bf Theorem 6.2.} {\it Let the group of automorphisms
$\tilde \alpha$ on the von Neumann factor $\cal M$ be obtained
through a Markovian cocycle perturbation of the Kolmogorov flow
generated by the quantum noise $j$ with the expectation defining
the probability distribution, which is invariant with respect to
the Markovian cocycle $w$.
Then there exists a subfactor
$\tilde {\cal M}\subset {\cal M}$ such that
the restriction $\tilde \alpha |_{\tilde {\cal M}}$
is the Kolmogorov flow generated by the quantum noise
$\tilde j$ which is isomorphic to the initial one.
The limit $\lim \limits _{t\to +\infty }w_{-t}=w_{-\infty }$
correctly defines a normal *-endomorphism $w_{-\infty }$
with the property $\tilde {\cal M}=w_{-\infty }({\cal M}),\
\tilde j_t=w_{-\infty }\circ j_t,\ t\in {\mathbb R}$.
}

Proof.

The first part of the theorem follows from Proposition 6.1.
Below we shall prove that there exists the limit
$\lim \limits _{t\to +\infty}
w_{-t}=w_{-\infty }$ defining the normal *-endomorphism on
${\cal M}$ with the properties we claimed.
Notice that $w_{-t-s}(y)=w_{-t}(y),$
$y\in {\cal M}_{t]},$ it follows that
the limit exists on the dense set of elements.
In the following we shall prove that it exists.
The formula $(U_t[x],[y])={\mathbb E}(\alpha _t(x)y^*),$ $
(W_t[x],[y])={\mathbb E}(w_t(x)y^*),\ x,y\in {\cal M},$
defines a strong continuous group
$U=(U_t)_{t\in \mathbb R}$ of unitary operators
in $\cal H$ and a unitary $U$-cocycle
$W=(W_t)_{t\in \mathbb R}$ correspondingly. Let ${\cal H}_{[t}$
be a subspace of $\cal H$ generated by the elements $[x],\
x\in {\cal M}_{[t}$. Then $W_t\xi =\xi ,\ \xi \in {\cal H}_{[t},\
t\geq 0,$ by means of the Markovian property for
$w$. Therefore,
$W$ is a Markovian $U$-cocycle in the sense of the definition of
Part 1 and the limit $s-\lim \limits _{t\to +\infty }W_{-t}=
W_{-\infty }$ exists due to Proposition 1.2.
Notice that $w_t(x)=W_txW_t^*,\ x\in {\cal M},\ t\in {\mathbb R}$.
It implies that $w_{-t-s}(x)\xi =W_{-t-s}xW_{-t-s}\xi=
W_{-t-s}xW_{-t}\xi \to W_{-\infty }xW_{-t}\xi ,\ x\in {\cal M},\
\xi \in {\cal H}_{[-t},\ s\to +\infty $. Hence,
the limit $\lim \limits _{t\to +\infty }w_{-t}(x)\xi $ exists
for a dense set of the elements $\xi \in {\cal H}$.
Because $w_{-t}$ is an automorphism, we get
$||w_{-t}(x)||=||x||$. Therefore, the strong limit
$s-\lim \limits _{t\to +\infty }w_{-t}(x)$ is defined
by the Banach-Steinhaus theorem for all
$x\in {\cal M}$. Moreover, the limiting map
$w_{-\infty }$ preserves the identity because all
$w_{-t}$ satisfy this property and
$||w_{-\infty }||=1$. In this way, the map
$w_{-\infty }$ is positive
(see \cite {Brat}). On the other hand,
$w_{-\infty }$ is a normal *-endomorphism for
it is a limit of the series of normal *-automorphisms
$w_{-t}$.
Thus, $w_{-\infty }$ is completely positive.
Notice that $w_{-\infty }({\cal M}_{0]})=
\tilde {\cal M}_{0]}$.
The Markovian property gives us
$w_{-t-s}(y)=w_{-t}(y),\ y\in {\cal M}_{[-t},\
s,t\geq 0$. Due to Proposition 6.1,
$\tilde j_{-t}=w_{-t}\circ j_{-t},\ t\geq 0$.
Hence, $\tilde j_{-t}=w_{-t-s}\circ j_{-t} =w_{-\infty }\circ
j_{-t},\ t\geq 0$. For positive values of time,
$\tilde j_t=j_t=w_{-\infty }\circ j_t,\ t\geq 0$,
by means of $w_{-t}(y)=w_{-\infty }(y)=y,\ t\geq 0,\
y\in {\cal M}_{[0}$. $\Box$

\section*{Acknowledgments} The author likes to thank Professor A.S.
Holevo for constant useful discussions resulted in
a significant improvement of the paper.

\end {document}